\documentclass[a4paper,12pt]{article}
\usepackage{amsthm,amsmath,amssymb,amsfonts}
\usepackage{epsfig}
\usepackage{graphicx,color}
\usepackage{pstricks}

\newtheorem{lem}{Lemma}
\newtheorem{cor}{Corollary}
\newtheorem{teo}{Theorem}

\newcommand{\vect}[1]{\mbox{\boldmath$#1$}}

\begin{document}

\title{Bijective mapping preserving intersecting antichains for $k$-valued cubes}
\author{Roman Glebov\\ 
\small{Universit\"{a}t Rostock, Institut f\"{u}r Mathematik, D-18051 Rostock, Germany}\\
\small{roman.glebov@uni-rostock.de}} 

\maketitle

\begin{abstract}
Generalizing a result of Miyakawa, Nozaki, Pogosyan and Rosenberg, 
we prove that there is a one-to-one correspondence 
between the set of intersecting antichains 
in a subset of the lower half of the $k$-valued $n$-cube 
and the set of intersecting antichains in the $k$-valued $(n-1)$-cube.
\end{abstract}

\section{Introduction}
Let $k$ and $n$ be positive integers with $k \geq 2$, and let $E = \{ 0, \ldots, k-1 \}$.
A \emph{$k$-valued $n$-cube} is the cartesian power $E^n$.
Let $ \vect{a} = (a_1,\ldots,a_n), \vect{b} = (b_1,\ldots,b_n) \in E^n$.
We write $\vect{a} \preceq \vect{b}$ if $a_i \leq b_i$ for all $i \in [ n ] := 
\{ 1, \ldots, n \}$.
We call $A \subseteq E^n$ an \textit{antichain} if
there are no different elements $ \vect{a}, \vect{b}$ of $A$
such that $\vect{a} \preceq \vect{b}$.
A family $A \subseteq E^n$ is \textit{intersecting} if for all
$\vect{a}, \vect{b} \in A$ there exists $i \in [ n ]$ such that $a_i + b_i
\geq k$.
This is a natural generalization of the binary case ($k=2$), 
where the elements of $E^n$ can be interpreted as the subsets of $[n]$
and an intersecting antichain is one consisting of pairwise intersecting sets.
The restriction in the definition applies also when $b=a$, 
so no $\vect{a} \in E^n$ with $a_i < \frac{k}{2}$ for all $i \in [n]$ 
is an element of any intersecting antichain, 
because then $a_i+a_i < k$ for all $i \in [n]$.

In the binary case, 
there is a bijective map from  the lower half of the $n$-cube 
onto the $(n-1)$-cube that preserves intersecting antichains in both directions \cite{r1}. 
Answering a question of Miyakawa \cite{r2}, 
we present a generalization to the $k$-valued case.
The proof is slightly simpler than that of \cite{r1} for the $k=2$ case.
More information on intersecting antichains can be found in \cite{r3}.

The {\it weight} of an element $\vect{a}\in E^n$, written $w(\vect{a})$, is defined by $w(\vect{a})=a_1+\ldots+a_n$.
For $0 \le t \le n(k-1)$ the $t$-th \textit{layer} $\mathcal{B}_t$ of $E^n$ is denoted by
$\mathcal{B}_t = \{ \vect{a} \in E^n : ~ w(\vect{a}) = t \} $. 
Now we define the "lower half" $L_n$ with restricted first entries.

Let $g = \lfloor \tfrac{n(k-1)}{2} \rfloor$ 
and notice that $g= \tfrac{1}{2}\left( nk-n-1 \right)$ if $n(k-1)$ is odd and 
$g= \tfrac{1}{2}n(k-1)$ otherwise. 
Let $C_i = \{ (a_1, \ldots , a_n) \in E^n : ~ a_1 = i \}$. 
Let
\[L_n = \begin{cases} (\mathcal{B}_{0} \cup \ldots \cup \mathcal{B}_{g}) \cap 
(C_0 \cup C_{k-1})
& \text{ if $n(k-1)$ is odd,}\\
((\mathcal{B}_{0} \cup \ldots \cup \mathcal{B}_{g-1}) \cap (C_0 \cup C_{k-1}))
\cup ( \mathcal{B}_{g} \cap C_0 )
& \text{ otherwise.}
\end{cases}\]

This set can be given also as follows: 
Let $g' = \lfloor \tfrac{n(k-1)-1}{2} \rfloor$,
and notice that $g' = \tfrac{1}{2}\left( nk-n-1 \right)=g$ if $n(k-1)$ is odd and 
$g' = \tfrac{1}{2}n(k-1)-1=g-1$ otherwise. 
Thus
\[ L_n =
\begin{cases}
(\mathcal{B}_{0} \cup \ldots \cup \mathcal{B}_{g'}) \cap (C_0 \cup C_{k-1})
& \mbox{ if $n(k-1)$ is odd,}\\
((\mathcal{B}_{0} \cup \ldots \cup \mathcal{B}_{g'}) \cap (C_0 \cup C_{k-1}))
\cup ( \mathcal{B}_{g'+ 1} \cap C_0 )
& \mbox{ otherwise.}
\end{cases}
\]
Hence, $g$ is the maximum weight of the elements of $L_n$ beginning with $0$. Similarly, 
$g'$ is the maximum weight of the elements of $L_n$ beginning with $k-1$.
Notice that $g + 1 + g' = n(k-1)$.

\section{A map from $L_n$ to $E^{n-1}$}

For $a \in E$, let $\overline{a} = k-1-a$. Define a map $ \varphi$ from
$L_n$ into $E^{n-1}$ by setting
\[
\varphi ((a_1,\ldots , a_n)) =
\begin{cases}
(a_2, \ldots, a_n)& \mbox{ if } a_1=0,\\
(\overline{a}_2,\ldots,\overline{a}_n)& \mbox{ if } a_1=k-1.
\end{cases}
\]
Obviously $\overline{\overline{a}} = a$ and
$a=b$ iff $\overline{a}=\overline{b}$. 
Concerning the weight $w$, note that
\[w(\varphi (\vect{a})) =
\begin{cases}
    w(\vect{a}) & \mbox{if } a_1 = 0,\\
    (k-1)(n-1)-(w(\vect{a})-(k-1)) & \mbox{if } a_1 = k-1.
\end{cases}\]

\begin{lem}
\label{weight}
For $\vect{a},\vect{b}\in L_n$ with $a_1 = 0$ and $b_1 = k-1$, 
we have \[w(\varphi (\vect{a})) < w(\varphi (\vect{b})).\]
\end{lem}

\proof

We have
\begin{multline*}
    w(\varphi (\vect{b}))
     = (k-1)(n-1)-(w(\vect{b})-(k-1)) 
     = n(k-1)-w(\vect{b}) \\
     = g+1+g' - w(\vect{b}) 
     \geq g+1 
     \geq w(\vect{a})+1 
     = w(\varphi (\vect{a}))+1 \\
     > w(\varphi (\vect{a})).
\end{multline*}
\qed

\begin{lem}
\label{inj} The map $\varphi$ is injective.
\end{lem}

\proof

Let $ \vect{a}, \vect{b} \in L_n, \vect{a} \not = \vect{b}$.
If $a_1 = b_1$, we obtain immediately from the definition of
$\varphi$ that $\varphi(\vect{a}) \not = \varphi(\vect{b})$.
If $a_1 \neq b_1$, 
w.l.o.g. $a_1 = 0$ and $b_1 = k-1$.
By Lemma \ref{weight}, $w(\varphi (\vect{b})) > w(\varphi (\vect{a}))$,
hence $\varphi(\vect{a}) \neq \varphi(\vect{b})$.  \qed

\begin{lem}
\label{surj}
The map $\varphi$ is surjective.
\end{lem}

\proof

We have to show that for all $\vect{b} = (b_1,\ldots,b_{n-1}) \in E^{n-1}$ there
exists an $\vect{a} \in L_n$
such that $\varphi (\vect{a}) = \vect{b}$.
We construct this $\vect{a}$ as follows: Let
\[\vect{a} =
\begin{cases}
    (0, b_1,\ldots,b_{n-1}) & \mbox{if } w(\vect{b}) \leq g,\\
    (k-1, \overline{ b}_1,\ldots,\overline{ b}_{n-1}) & \mbox{if } w(\vect{b}) > g.
\end{cases}\]
If $w(\vect{b}) \leq g$, then $w(\vect{a}) =
w(\vect{b}) \leq g$.
If $w(\vect{b}) > g$, then
$w(\vect{a}) = k-1 + ((k-1)(n-1)- w(\vect{b})) < n(k-1)-g = g'+1 $, hence
$w(\vect{a}) \leq g'$. Thus in both cases $\vect{a} \in
L_n$, and obviously $\varphi (\vect{a}) = \vect{b}$. \qed

\begin{cor}
\label{bi}
\textit{The map $\varphi:L_n\to E^{n-1}$ is a bijection.}
\end{cor}

\begin{lem}
\label{pre}
Both $\varphi$ and its inverse preserve intersecting antichains.
\end{lem}

\proof

Due to the definition of an intersecting antichain,
it is sufficient to prove the lemma for antichains $A$ with $|A| \in \{ 1,2\}$.

Let
$ \vect{a}, \vect{b} \in L_n$ and let
$\{ \vect{a}, \vect{b} \}$ be an intersecting antichain.

If $a_1 = b_1 = 0$, then obviously $\{\varphi (\vect{a}), \varphi (\vect{b}) \}$
is an intersecting antichain.

If $a_1 = b_1 = k-1$, then
\begin{align*}
    w(\varphi (\vect{a})) + w(\varphi (\vect{b}))
    & = (k-1)(n-1)-(w(\vect{a})-(k-1)) \\
    & \qquad + (k-1)(n-1)-(w(\vect{b})-(k-1))\\
    & \geq  2n(k-1)-2\left\lfloor \frac{n(k-1)-1}{2} \right\rfloor \\
    & >  (k-1)(n-1).
\end{align*}
Thus, there exists $i \in \{2,\ldots, n \} $ such that
$\overline{ a}_i + \overline{ b}_i \geq k$, and
hence $\{\varphi (\vect{a}), \varphi (\vect{b}) \}$ is intersecting.
Furthermore, if $ \vect{a} = \vect{b}$, obviously 
$\{\varphi (\vect{a}), \varphi (\vect{b}) \} = \{\varphi (\vect{a}) \}$ is an antichain.
If $ \vect{a} \neq \vect{b}$, by the antichain property, 
there are $i,j \in \{2,\ldots, n \} $ with $a_i < b_i$ and $a_j > b_j$.
Thus $\overline{ a}_i > \overline{ b}_i$ and $\overline{ a}_j < \overline{ b}_j$, and hence
$\{\varphi (\vect{a}), \varphi (\vect{b}) \}$
is an antichain.

If $a_1 \not = b_1$, then we may assume $a_1 = 0$ and $b_1 = k-1$. 
Obviously $ \vect{a} \neq \vect{b}$.
By Lemma \ref{weight}, $w(\varphi (\vect{a})) < w(\varphi (\vect{b}))$, and thus
$\varphi (\vect{a}) \not \succeq  \varphi (\vect{b})$.
Since $\{ \vect{a}, \vect{b} \}$ is intersecting,
there existsn $i \in \{ 2,\ldots, n \} $, such that $a_i + b_i \geq k$.
Thus $ \overline{ b}_i = k-1 - b_i < a_i$, hence
$\varphi(\vect{a}) \not \preceq  \varphi (\vect{b})$.
Consequently $\{\varphi (\vect{a}), \varphi (\vect{b}) \}$ is an antichain.
Since $\{ \vect{a}, \vect{b} \}$ is an antichain,
there existsn $i \in \{ 2, \ldots, n \} $, such that $a_i > b_i$,
so $ a_i + \overline{ b}_i = a_i + k-1 - b_i > k-1$, and hence
$\{\varphi (\vect{a}), \varphi (\vect{b}) \}$ is
intersecting.

Now let
$ \vect{a}, \vect{b} \in E^{n-1}$ and let
$\{ \vect{a}, \vect{b} \}$ be an intersecting antichain.
By the proof of Lemma \ref{surj}, for
$\vect{b} \in E^{n-1}$,
\[\varphi^{-1}(\vect{b}) =
\begin{cases}
    (0, b_1,\ldots,b_{n-1}) & \mbox{if } w(\vect{b}) \leq g,\\
    (k-1, \overline{ b}_1,\ldots,\overline{ b}_{n-1}) & \mbox{if } w(\vect{b}) > g.
\end{cases}\]

If $w(\vect{a}) \leq g$ and $w(\vect{b}) \leq g$,
then obviously $\{\varphi^{-1} (\vect{a}), \varphi^{-1} (\vect{b}) \}$
is an intersecting antichain.

If $w(\vect{a}) > g$ and $w(\vect{b}) > g$, then
the first entry of both
$\varphi^{-1} (\vect{a})$ and $\varphi^{-1} (\vect{b}) $ is $k-1$,
so $\{\varphi^{-1} (\vect{a}), \varphi^{-1} (\vect{b}) \}$
is intersecting.
Furthermore, if $ \vect{a} = \vect{b}$, obviously 
$\{\varphi^{-1} (\vect{a}), \varphi^{-1} (\vect{b}) \} = 
\{\varphi^{-1} (\vect{a}) \}$ is an antichain.
If $ \vect{a} \neq \vect{b}$, there are $i,j \in [n-1]$ with $a_i < b_i$ and $a_j > b_j$,
thus $\overline{ a}_i > \overline{ b}_i, \overline{ a}_j < \overline{ b}_j$, and hence
$\{\varphi^{-1} (\vect{a}), \varphi^{-1} (\vect{b}) \}$
is an antichain.

In the remaining case, we may assume $w(\vect{a}) \leq g$ and $w(\vect{b}) > g$.
Obviously $ \vect{a} \neq \vect{b}$.
The first entry of $\varphi^{-1} (\vect{a})$ is $0$
and the first entry of $\varphi^{-1} (\vect{b})$ is $k-1$,
so $\varphi^{-1} (\vect{a}) \not \succeq  \varphi^{-1} (\vect{b})$.
Since $\{ \vect{a}, \vect{b} \}$ is intersecting,
there existsn $i \in [n-1]$ such that $a_i+b_i \geq k$. 
Thus
$a_i \geq k-b_i = \overline{b}_i +1 > \overline{b}_i$, 
and hence
$\varphi^{-1} (\vect{a}) \not \preceq  \varphi^{-1} (\vect{b})$.
Consequently $\{\varphi^{-1} (\vect{a}), \varphi^{-1} (\vect{b}) \}$ is an antichain.
Since $\{ \vect{a}, \vect{b} \}$ is an antichain,
there existsn $i \in [n-1]$ such that $a_i > b_i$,
thus $ a_i + \overline{ b}_i = a_i + k-1 - b_i > k-1$, hence
$\{\varphi^{-1} (\vect{a}), \varphi^{-1} (\vect{b}) \}$ is
intersecting.  \qed

From Corollary \ref{bi} and Lemma \ref{pre} 
we immediately obtain the main result of this note.

\begin{teo}
\label{satz}
\textit{The map $\varphi$ is bijective
and preserves intersecting antichains in both directions.}
\end{teo}

\section{Maiximum Size of an Antichain and an Intersecting Antichain in $E^n$ and $L_n$}

To show an application of Theorem~\ref{satz}, 
we first estimate the maximum size of an intersecting antichain in $E^n$.

\begin{teo}
\label{iaen}
\textit{The map $\varphi$ is bijective
and preserves intersecting antichains in both directions.}
\end{teo}

\proof 
Let $W$ be a maximum intersecting antichain of size $m$.
Set 
\[s:=\min \{t: W \cap \mathcal{B}_t \neq \emptyset\},\]
\[W':=\left( W^{\succeq }\cap\mathcal{B}_{s+1}\right)\cup\left( W \backslash \mathcal{B}_{s}\right).\]
A direkt check shows that $W'$ is an intersecting antichain.

\section{Remarks}

In the definition of $L_n$, we can replace $C_0$ by $C_i$ and
$C_{k-1}$ by $C_{k-1-i}$ with $0 \leq i < \frac{k-1}{2} $. We obtain
\[
L_{n,i} =
\begin{cases}
(\mathcal{B}_{0} \cup \ldots \cup \mathcal{B}_{g}) \cap (C_i \cup C_{k-1-i})
& \text{ if } n(k-1) \text{ is odd,}\\
((\mathcal{B}_{0} \cup \ldots \cup \mathcal{B}_{g-1}) \cap (C_i \cup C_{k-1-i}))
\cup ( \mathcal{B}_{g} \cap C_i )
& \text{ otherwise.}
\end{cases}
\]
The analogue on $L_{n,i}$ of the map $\varphi$ on $L_n$ also is a bijection to $E^{n-1} $
and preserves intersecting antichains.
The only place where the proof is not completely identical
is the case $a_1 = b_1 = k-1-i$ in the first direction of Lemma \ref{pre}.
In this case, we have
\begin{align*}
    w(\varphi (\vect{a})) + w(\varphi (\vect{b}))
    & = (k-1)(n-1)-(w(\vect{a})-(k-1-i) \\
    & \qquad + (k-1)(n-1)-(w(\vect{b})-(k-1-i))\\
    & \geq  2(n-1)(k-1) - 2\left\lfloor \frac{n(k-1)-1}{2} \right\rfloor +
    2(k-1-i) \\
    & >  2(k-1)(n-1) - n(k-1) + (k-1)\\
    & =  (k-1)(n-1).
\end{align*}

Furthermore, $g$ can be replaced by $g+z$ and $g'$ by $g'-z$
with $z \in \{ 0, \ldots, g' \} $, such that
\[
L_{n}^z :=
((\mathcal{B}_{0} \cup \ldots \cup \mathcal{B}_{g+z}) \cap C_0 ) \cup
((\mathcal{B}_{0} \cup \ldots \cup \mathcal{B}_{g'-z}) \cap C_{k-1}).
\]
As in the definition in Lemma \ref{surj}, for
$\vect{b} \in E^{n-1}$ we have
\[\varphi^{-1}(\vect{b}) =
\begin{cases}
    (0, b_1,\ldots,b_{n-1}) & \mbox{if } w(\vect{b}) \leq g+z,\\
    (k-1, \overline{ b}_1,\ldots,\overline{ b}_{n-1}) & \mbox{if } w(\vect{b}) > g+z.
\end{cases}\]

\bigskip

{\bf Acknowledgement.} We are grateful to Konrad Engel and Florian
Pfender for helpful suggestions. We are also thankful to Thomas Kalinowski
and Antje Kiesel for carefully reading the paper and proposing many correc-
tions greatly improving both its English and general readability.

\end{document}